\documentclass[12pt,a4paper, reqno]{amsart}
\usepackage[utf8]{inputenc}	
\usepackage[T1]{fontenc}
\usepackage[in,headings]{fullpage}
\usepackage{amsbsy, amscd, amsfonts, amsmath, amsrefs, amssymb, amsthm}
\usepackage{graphicx}
\usepackage{float}
\usepackage{color}   %May be necessary if you want to color links
\usepackage[colorlinks=true,linkcolor=blue,citecolor=blue,urlcolor=blue]{hyperref}
\usepackage{comment}
\excludecomment{A}
\usepackage{epigraph}

\newtheorem{theorem}{Theorem}%[section]

\newtheorem{corollary}[theorem]{Corollary}
\newtheorem{proposition}[theorem]{Proposition}
\newtheorem{lemma}[theorem]{Lemma}

\theoremstyle{definition}
\newtheorem{definition}[theorem]{Definition}
\newtheorem{remark}[theorem]{Remark}

\theoremstyle{remark}

\newcommand{\C}{\mathbf{C}}
\newcommand{\Z}{\mathbf{Z}}

\newcommand{\R}{\mathbf{R}}
\newcommand{\N}{\mathbf{N}}

\renewcommand{\Re}{\mathop{\mathrm{Re}}\nolimits}
\renewcommand{\Im}{\mathop{\mathrm{Im}}\nolimits}

\newfont{\cmbsy}{cmbsy10}
\newfont{\cmmib}{cmmib10}
\newcommand{\Orden}{\mathop{\hbox{\cmbsy O}}\nolimits}

\DeclareMathOperator*{\card}{card} 
\newcommand{\Res}{\mathop{\operatorname{Res}}}

\begin{document}

\title{Asymptotic properties of zeros of Riemann zeta function}
\author[Arias de Reyna]{Juan Arias de Reyna}
\address{%
Universidad de Sevilla \\ 
Facultad de Matem\'aticas \\ 
c/Tarfia, sn \\ 
41012-Sevilla \\ 
Spain.} 
\email{arias@us.es, ariasdereyna1947@gmail.com}

\author[Meyer]{Yves Meyer}
\address{%
13, rue des Potiers \\ 
92260 Fontenay-aux-Roses \\ 
France.} 
\email{Yves.Meyer305@orange.fr}

%\epigraph{Papers should include more side remarks, open questions, and such. Very often, these are more interesting than the theorems actually proved.\\
%\textsc{Jean-Pierre Serre}}

% AMS subject classifications (used in AMS journals)

% AMS keywords (used in AMS journals)

%\date{\today, \texttt{104-Zeros-v12.tex}}

\begin{abstract}
We try to define the sequence of zeros of the Riemann zeta function by an intrinsic property.  Let  $(z_k)_{k\in \mathbb{N}}$ be the  sequence of nontrivial zeros of $\zeta(s)$ with positive imaginary part. We write $z_k= 1/2+i\tau_k$ (RH says that these $\tau_k$ are all  real). Then the sequence  $(\tau_k)_{k\in \mathbb{N}},$ satisfies the following asymptotic relation
\[\sum_{k\in\mathbb{N}}\frac{2x}{x^2+\tau_k^2}\simeq \frac12\log\frac{x}{2\pi}+\sum_{n=1}^\infty \frac{a_n}{x^n},\,\,x\to +\infty\]
 where $a_{2n+1}=2^{-2n-2}(8-E_{2n})$, $a_{2n}=(1-2^{-2n+1})B_{2n}/(4n).$ Are there other sequences $(\alpha_k)_{k\in \mathbb{N}},$ of real or complex numbers enjoying this property? These problems are addressed in this note. 
\end{abstract}

\maketitle

\section{Introduction} 
The goal of this paper is to characterize the sequence of zeros of zeta as a whole by its properties. We do not assume the Riemann hypothesis. 

 As usual $ s=\sigma+i t,$ $\xi(s)=\frac{s(s-1)}{2}\Gamma(s/2)\pi^{-s/2}\zeta(s)$ and $\Xi(s)=\xi(1/2+is).$  
Then $\Xi(s)$ is an even entire function, $\Xi(\tau_n)=0,$  and the product formula for the zeta function can be written
\[\zeta(\tfrac12+t)=\frac{\Xi(0)\pi^{\frac14+\frac{t}{2}}}{\frac12(t^2-\frac14)\Gamma(\frac14+\frac{t}{2})}\cdot\prod_{\substack{n\in \mathbb{N}}}\Bigl(1+\frac{t^2}{\tau_n^2}\Bigr).\]
For $t\to+\infty$ the left hand side converges to $1$ exponentially. Both $\Gamma(\frac14+\frac{t}{2})$ and the infinite product $\prod_{\substack{n\in \mathbb{N}}}\Bigl(1+\frac{t^2}{\tau_n^2}\Bigr)$  tend to $\infty$ as $t$ tends to $+\infty.$ Therefore  the two factors on the right hand side  have to be delicately balanced. This is the core of the properties of the zeros we want to use.

This approach is inspired by an equation 
\[\sum_{n\in\mathbb{N}}\log\Bigl(1+\frac{t^2}{\tau_n^2}\Bigr)=\Bigl(\frac{t}{2}+\frac74\Bigr)\log\frac{t}{2\pi}-\frac{t}{2}+2\log\pi+\frac32\log2-\log\Xi(0)+O(t^{-1})\]
that we find in Riemann's Nachlass.\footnote{(Sheet 77 of Number Theory). We have changed slightly Riemann's notations.}

In Section \ref{suno} we select some of the properties satisfied by the sequence of zeros of $\Xi(t)$. Then (Section \ref{inverse}) we introduce the concept of Riemann sequences, those that satisfy some of the properties of  zeros of $\Xi(t)$ with positive real part.

The zeros of $\Xi(t)$ give the first example of a Riemann sequence. In Section \ref{S:newzeta} (see Theorem \ref{T:mainconst}) we construct an infinity of examples of Riemann sequences. It is expected all these examples are \textbf{complex} Riemann sequences. In Section \ref{S:numeric} we construct a concrete example of a complex Riemann sequence. 
The construction of these new Riemann sequences follows from the results of \cite{M}. To this end we reduce the problem to the construction of certain crystalline measures. This implies the existence  of some meromorphic  functions that share with zeta the functional equation, the pole at $s=1$ as a unique pole, that are the sum of a convergent Dirichlet series for $\sigma>\sigma_0$ and converge to $1$ exponentially  for $s\to+\infty$. We do not expect these functions to have an Euler product.

In Section \ref{sproperties} we show many properties that are common to all Riemann sequences. Of course they are known for the zeta function. See the works of Chakravarty \cite{C1},  \cite{C3} and \cite{C4},  and Voros \cite{V}, and \cite{V2}.
The point is that we obtain all these properties for all Riemann sequences.

The Riemann hypothesis yields a  real  Riemann sequence $(\tau_n)$. Another real Riemann sequence will imply the existence of a meromorphic function with the same functional equation as the $\zeta$ function and zeros on the critical line. Therefore, we expect that $(\tau_n)$ is the only \emph{real} Riemann sequence. This is the characterization we were looking for.

We have several interesting questions about Riemann sequences:
\begin{itemize}
\item[(1)] Are there some ``real'' Riemann sequences? 
\item[(2)] Is there only one real Riemann sequence?
\item[(3)] Is any Riemann sequence associated to a Dirichlet series?
\end{itemize}

One nice properties of the first two questions is that they do not imply RH and  positive answers are believed in by believers.  Since RH is such a difficult problem it is good to have  some different and related problems (one of Polya's advice).

\subsection{Notations}
We make the convention throughout that $\tfrac12+i\tau$ shall run through the complex zeros of the function $\zeta(s)$, the Riemann's zeta-function. Thus some of the values of $\tau$ will be complex if the Riemann hypothesis is false.  We form a sequence $(\tau_n)$ with all the $\tau$ with $\Re\tau>0$ ordered  in such a way that $0<\Re\tau_1\le \Re\tau_2\le \cdots$ Furthermore, if multiple zeros of $\zeta(s)$ exist, we repeat the corresponding terms  in the sequence according to their  order of multiplicity. 

$C_0=0.577216\dots$ denotes the Euler constant. As usual  
$E_n$ and $B_n$ denotes the ordinary Euler and Bernoulli numbers
\begin{equation}
\frac{2e^t}{e^{2t}+1}=\sum_{n=0}^\infty E_n\frac{t^n}{n!},\qquad
\frac{t}{e^t-1}=\sum_{n=0}^\infty B_n\frac{t^n}{n!}.
\end{equation}

Given a function $f(z)$ defined on a half-plane $\Re z>a$ we shall write
\begin{equation}
f(z)\simeq\sum_{n=0}^\infty \frac{a_n}{z^n}
\end{equation}
if for each natural number $N$ and any $0<\varepsilon<\frac{\pi}{2}$ we have
\begin{equation}\label{defasym}
\lim_{z\to\infty}\Bigl(f(z)-\sum_{n=0}^N \frac{a_n}{z^n}\Bigr)z^N=0,\qquad |\arg z|\le
\frac{\pi}{2}-\varepsilon.
\end{equation}
Also we shall put $f(z)\simeq g(z)+\sum_n\frac{a_n}{z^n}$ if $f(z)-g(z)\simeq
\sum_n\frac{a_n}{z^n}$. 
We shall use $\sim$ to denote another type of asymptotic relationship.

\section{Some properties  of Riemann zeta function zeros.}\label{suno}

Our first theorem is a precise version of a Riemann equation we found in Riemann's Nachlass. 
\begin{theorem}\label{one}
We have the following asymptotic expansion in the sense defined in   \eqref{defasym}
\begin{equation}\label{logexpan}
\sum_{n=1}^\infty\log\Bigl(1+\frac{z^2}{\tau_n^2}\Bigr)\simeq
\frac{z}{2}\log\frac{z}{2\pi}-\frac{z}{2}+\frac{7}{4}\log z+A-
\sum_{n=1}^{\infty}\frac{a_{n+1}}{n z^n},\qquad z\to\infty
\end{equation}
where $A=\frac14\log\frac{\pi}{2}-\log\Xi(0)$ and the coefficients are given by
\begin{equation}\label{coeff}
a_{2n+1}=\frac{8-E_{2n}}{2^{2n+2}},\qquad a_{2n}=\Bigl(1-\frac{1}{2^{2n-1}}
\Bigr)\frac{B_{2n}}{4n}.
\end{equation}
\end{theorem}

\begin{proof}

The product formula  $\Xi(t)=\Xi(0)\prod_{n=1}^\infty\left(1-\frac{t^2}{\tau_n^2}\right)$ may be written as
\begin{equation}\label{productzeta}
\zeta(s)=\Xi(0)\frac{\pi^{\frac{s}{2}}}{(s-1)\Gamma(1+s/2)}\prod_{n=1}^\infty
\Bigl(1+\frac{(s-1/2)^2}{\tau_n^2}\Bigr).
\end{equation}\label{E:log}
Put $s=\tfrac12+z$ in  \eqref{productzeta} and take log on both sides to get
\begin{multline}\label{logonall}
\log\zeta(\tfrac12+z)=\\
\sum_{n=1}^\infty\log\Bigl(1+\frac{z^2}{\tau_n^2}\Bigr)+
\Bigl(\frac{z}{2}+\frac14\Bigr)
\log\pi-\log\Bigl(z-\frac12\Bigr)-\log\Gamma\Bigl(\frac{z}{2}+\frac{5}{4}\Bigr)+\log\Xi(0)
\end{multline}
Here the equality is true if for  $\Re z>a>1/2$ we use the usual analytic branch of $\log\zeta(\frac12+z)$ and $\log\Gamma(\frac{z}{2}+\frac54).$  
We define  $\log(z-\frac12)$ by choosing its  argument $\le \pi/2$. For the series $\sum_{n=1}^\infty\log(1+\frac{z^2}{\tau_n^2})$ we  notice that $1+z^2/\tau_n^2$ does not vanish for $\Re z>1/2$ or $|\Im z|<1/2$. We can take the branch of $\log(1+\frac{z^2}{\tau_n^2})$ that vanish at $z=0$, then the series converge since $\tau_n\to\infty$ and $\sum|\tau_n|^{-2}<+\infty$. With these definitions of logarithms the equality \eqref{logonall} holds.

From \cite{L}*{p.~32} we quote the well known expansion
\begin{equation}\label{eq7}
\log\Gamma(z+a)\simeq(z+a-\tfrac12)\log z-z+\tfrac12\log2\pi+\sum_{n=1}^\infty
(-1)^{n+1}\frac{B_{n+1}(a)}{n(n+1)}\frac{1}{z^n}.
\end{equation}
By the Taylor expansion 
\begin{equation}\label{eq8}
\log\Bigl(z-\frac12\Bigr)=\log z+\log\Bigl(1-\frac{1}{2z}\Bigr)\simeq\log z-\sum_{n=1}^\infty\frac{1}{n\,2^n}\frac{1}{z^n}.
\end{equation}
Finally $\log\zeta(\tfrac12+z)=\Orden(2^{-\Re z})$ which implies that  in 
our sense $\log\zeta(\tfrac12+z)\simeq 0$. This is why we need $\varepsilon>0$ in \eqref{defasym}.

Then from \eqref{logonall} we get
\begin{multline}
\sum_{n=1}^\infty\log\Bigl(1+\frac{z^2}{\tau_n^2}\Bigr)\simeq
\frac{z}{2}\log\frac{z}{2\pi}-\frac{z}{2}+\frac74\log z+\frac14\log \frac{\pi}{2}-\log\Xi(0)+\\
+\sum_{n=1}^\infty \Bigl(-\frac{1}{n\,2^n}
+(-1)^{n+1}\frac{2^nB_{n+1}(5/4)}{n(n+1)}\Bigr)\frac{1}{z^n}
\end{multline}

Hence the first few  terms of the expansion \eqref{logexpan} are 
\begin{multline}
\sum_{n=1}^\infty \Bigl(-\frac{1}{n\,2^n}
+(-1)^{n+1}\frac{2^nB_{n+1}(5/4)}{n(n+1)}\Bigr)\frac{1}{z^n}=\\
=
-\frac{1}{48\,z}-\frac{9}{32\, z^2}+\frac{7}{5760\, z^3}-\frac{3}{256\, z^4}-\frac{31}
{80640\, z^5}-\frac{23}{512\, z^6}+\cdots
\end{multline}
The values of the coefficients \eqref{coeff} follows from this expression 
and the next Lemma.
\end{proof}

\begin{lemma}
For $n\ge0$ we have the following values of the Bernoulli polynomials at $x=5/4$
\begin{equation}
B_{2n+1}(5/4)=\frac{2n+1}{2^{4n+2}}(4-E_{2n}), \quad
B_{2n}(5/4)=\frac{8n-(2^{2n}-2)B_{2n}}{2^{4n}},
\end{equation}
\end{lemma}
\begin{proof}
By definition
\begin{equation}
\frac{t e^{xt}}{e^t-1}=\sum_{n=0}^\infty B_n(x)\frac{t^n}{n!}.
\end{equation}
and our Lemma follows from the following elementary equalities
\begin{align*}
\frac{te^{5t/4}}{e^t-1}&=te^{t/4}+\frac{t}{4\sinh\frac{t}{4}}-\frac{t}{4\cosh\frac{t}{4}}=\\
&=\sum_{n=0}^\infty \frac{t^{n+1}}{4^n n!}+ \sum_{n=0}^\infty (2-2^{2n})\frac{B_{2n}t^{2n}}{4^{2n} (2n)! }-\sum_{n=0}^\infty
\frac{E_{2n} t^{2n+1}}{4^{2n+1} (2n)!}.\qedhere
\end{align*}
\end{proof}

\begin{theorem}\label{two}
If the nontrivial zeros of the zeta function with positive imaginary part  are written $z_k=\tfrac12+i\tau_k, \tau_k>0,$ then the sequence $(\tau_k)_{k\in \mathbb{N}},$ satisfies  the following asymptotic expansion
in the sense of \eqref{defasym}
\begin{equation}\label{main}
\sum_{k\in\mathbb{N}}\frac{2z}{z^2+\tau_k^2}\simeq\frac12\log\frac{z}{2\pi}+\sum_{n=1}
^\infty\frac{a_n}{z^n} 
\end{equation}
where the coefficients, even $a_1=\frac74$ are given by \eqref{coeff}.
\end{theorem}

\begin{proof}
It is well known that a relation such as \eqref{defasym} $f(z)$ being holomorphic 
may be differentiated. 
Therefore our Theorem follows by differentiating the result of  Theorem \ref{one}.
We must check also that the coefficient $a_1=\frac74$ is given by \eqref{coeff}.
\end{proof} 

The first few terms of  expansion \eqref{main} are
\begin{equation}
\sum_{n\in\mathbb{N}}\frac{2z}{z^2+\tau_n^2}\simeq\frac12\log\frac{z}{2\pi}
+\frac{7}{4\,z}+\frac{1}{48\,z^2}+\frac{9}{16\,z^3}-\frac{7}{1920\,z^4}+\frac{3}{64\,z^5}+\cdots
\end{equation}

In the rest of the paper we study sequences of complex numbers sharing some properties with the $(\tau_n)$.

\section{Riemann's Sequences.}\label{inverse}

In Section \ref{suno} we have seen that the zeros of Riemann zeta function have property \eqref{main}. We define a Riemann sequence as a sequence of complex number sharing some properties with the zeros of $\Xi(t)$. We will see that there are many examples of Riemann sequences. But we still believe that the zeros of $\Xi(t)$ is the unique \emph{real} Riemann sequence.

Given a sequence $(\alpha_n)$ of complex numbers  and a complex number $\alpha$ we define the multiplicity of $\alpha$ as $\nu(\alpha)=\card\{n\in\N : \alpha_n=\alpha\}$. 

\begin{definition}\label{D:RS}
A sequence of complex number $(\alpha_n)$ is a Riemann's sequence if:

\begin{enumerate}

\item[(a)] $\displaystyle{0<1<\Re\alpha_1\le \Re\alpha_2\le\Re\alpha_3\le \cdots}$.

\item[(b)] The multiplicity of each complex number is finite and if $\alpha$ is not 
real $\nu(\alpha)=\nu(\overline{\alpha})$.

\item[(c)] There is a constant $C$ such that  $|\Im\alpha_n|<C$ for each value of $n\in\N$.

\item[(d)] We have $|\arg(\alpha_n)|<\pi/4$, or equivalently $|\Im(\alpha_n)|<\Re(\alpha_n)$.  

\item[(e)] The sequence satisfies 
\begin{equation}\label{mainalpha}
\sum_{k=1}^\infty\frac{2z}{z^2+\alpha_k^2}\simeq\frac12\log\frac{z}{2\pi}+\sum_{n=1}
^\infty\frac{a_n}{z^n} 
\end{equation}
\noindent where the $a_n$ are the numbers defined in \eqref{coeff}.
\end{enumerate}
\end{definition}

Of course $(\tau_n)$ is a Riemann sequence. The restriction $1<\Re\alpha_1$ is only a technical hypothesis. It simplifies some  of the following reasoning and is satisfied by the sequence $(\tau_n)$ related to the zeros of zeta.

Let $(\alpha_k)_{k\in \mathbb{N}}$ be a real Riemann's sequence. If $z=1/\epsilon$ and $\epsilon \to 0$ we obtain the asymptotic expansion 
\begin{equation}
\sum_{k\in\mathbb{N}}\frac{2\epsilon}{1+\epsilon^2\alpha_k^2}\simeq\frac12\log\frac{1/\epsilon}{2\pi}+\sum_{n=1}
^\infty a_n{\epsilon}^n 
\end{equation}
which can be viewed as a particular   example of a more general problem: what is the behavior of $\sum_{k\in\mathbb{N}} \epsilon \phi(\epsilon \alpha_k)$ as $\epsilon \to 0,  \phi$ being a non negative function which is $\mathcal{O}(x^{-2})$ at infinity?  An example is $\phi=\chi_{[0, 1]}$ and we then count the number of $\alpha_k$ in $[0, 1/\epsilon].$

\section{Crystalline measures}\label{sconst}

A measure $\mu$ on $\R^n$ is a \emph{crystalline measure} if it is a tempered distribution and both $\mu$ and its Fourier transform  have a locally finite support. The main example of a crystalline measure in $\R$ is the Dirac comb $\mu=\sum_{n\in\Z}\delta_n$ whose Fourier transform $\widehat\mu$ is the same $\widehat\mu=\mu$. By translating an rescaling this example it is easy to see that for any real numbers $a$ and $b$ with $b\ne0$, we have $\widehat\mu=\nu$ for the pair of measures
\begin{equation}\label{DiracCombs}
\mu:=\sum_{n\in\Z}\delta_{a+bn},\qquad \nu:=\frac1{|b|}\sum_{m\in\Z}e^{-2\pi i a m/b}\delta_{m/b}.
\end{equation}
Linear combinations of these examples do not exhaust the vector space of crystalline measures.

There is a connection between crystalline measure and zeta functions with a functional equation. This is due to Hans Hamburger \cite{H}. 
\begin{theorem}\label{T:cryszeta}
Let $\lambda_k$ for $k\in\Z$ be an odd sequence  $\lambda_{-k}=-\lambda_{k}$ with $1\le\lambda_1<\lambda_2<\cdots$ with $\lambda_k\to+\infty$ as $k\to+\infty$ and $c_k$ with $k\in \Z$ an even sequence of non null complex numbers. Assume also that $\sum_{k=1}^\infty |c_k|\lambda_k^{-\sigma_0}<+\infty$ for some $\sigma_0>0$. Let 
\[\mu:=\sum_{k\in \Z}c_k\delta_{\lambda_k},\qquad F(s):=\sum_{n=1}^\infty \frac{c_n}{\lambda_n^s}.\]
Let us observe that $\mu$ is a tempered measure and $F(s)$ is analytic for $\Re s>\sigma_0$. Then the following assertions are equivalent:
\begin{itemize}
\item[(a)] The measure $\mu$ is crystalline and $\widehat\mu=\mu$.
\item[(b)] The function $F(s)$ extends to a meromorphic function in the complex plane, with a unique pole at $1$ with residue $c_0$. The function $(s-1)F(s)$ is entire of order 1 and $F(s)$ satisfies the functional equation
\[\pi^{-\frac{s}{2}}\Gamma(\tfrac{s}{2})F(s)= \pi^{-\frac{1-s}{2}}\Gamma(\tfrac{1-s}{2})F(1-s).\]
\end{itemize}
\end{theorem}

We will use Theorem \ref{T:cryszeta} to construct  some special zeta functions sharing many properties with $\zeta(s)$. The main step is the construction of the corresponding  crystalline measures. Our construction follows the line of the one in \cite{M}*{Lemma 7}, where a crystalline measure is obtained such that $\widehat\mu=\mu$ and its restriction to $[-T,T]$ is equal to $0$. We modify here this construction to get some additional properties (see Proposition \ref{P:250131-4}) and to get concrete numerical examples.

For $N\in\N$ let  $\mathcal D_N$ be the space of measures
\[\mu=\sum_{n\in\Z}c_n\delta_{n/N},\qquad c_{n+N^2}=c_N,\quad c_n\in\C.\]
\begin{proposition}\label{P-250131-2}
Each element in $\mathcal D_N$ is a tempered distribution. The Fourier transform acts on this space $\mathcal F\colon\mathcal D_N\to\mathcal D_N$. In fact
\[\mathcal F(\mu)=\sum_{m\in\Z}P(m)\delta_{m/N},\]
where $P(x)$ is the trigonometric polynomial in \eqref{E:250130-1}.
\end{proposition}
\begin{proof}
Each integer $n$ can be written as $n=mN^2+r$ with $0\le r\le N^2-1$. Therefore,
\[\mu=\sum_{n\in\Z}c_n\delta_{n/N}=\sum_{r=0}^{N^2-1} c_r\sum_{m\in\Z}\delta_{mN+r/N}.\]
It follows that 
\[\widehat\mu=\sum_{r=0}^{N^2-1} c_r\mathcal F\Bigl(\sum_{m\in\Z}\delta_{mN+r/N}\Bigr)=
\sum_{r=0}^{N^2-1} c_r\frac1N\sum_{m\in\Z} e^{-2\pi i \frac{mr}{N^2}}\delta_{m/N}=\sum_{m\in\Z}P(m)\delta_{m/N},\]
where 
\begin{equation}\label{E:250130-1}
P(x)=\sum_{r=0}^{N^2-1} \frac{c_r}{N} e^{-2\pi i\frac{r}{N^2}x}.
\end{equation}
and obviously we have $P(x)=P(x+N^2)$. 
\end{proof}
The space $\mathcal D_N$ have dimension $N^2$ and the Fourier transform is given by 
\[P(m)=\frac1N\sum_{m=0}^{N^2-1}c_n e^{-2\pi i\frac{mn}{N^2}}, \qquad 0\le n\le N^2-1.\]
This is just the finite Fourier transform on $\Z/N^2\Z$. In general the finite Fourier transform is defined on $E_M$, the space of complex functions $f\colon\Z/M\Z\to\C$, then  $F\colon E_M\to E_M$ the finite Fourier transform is defined by 
\[F(f)=\widehat f, \qquad \text{where}\quad \widehat f(n)=\frac{1}{\sqrt{M}}\sum_{m\in\Z/M\Z}f(m)e^{-2\pi i \frac{mn}{M}}.\]
Our transformation corresponds to the case $M=N^2$.

Then it is known that $F$ is an isometry in $E_M$ endowed with the norm $\Vert f\Vert_2=\sum_{n\in\Z/M\Z} |f(n)|^2$. We have also $F^2(f)=g$ with $g(n)=f(-n)$, and $F^4=I$. Therefore the eigenvalues of $F$ satisfies $\lambda^4=1$. Let $V_\lambda$ the subspace of $E_M$ with $F(f)=\lambda f.$ Then the dimension of these spaces equals also to the multiplicities of the corresponding eigenvalue and is given by the table taken from \cite{Mc}

\begin{table}[htp]
%\caption{default}
\begin{center}
\begin{tabular}{lllll}
\hline
\\
M & $\dim(V_1)$ & $\dim(V_{-1})$ & $\dim(V_{-i})$ & $\dim(V_{i})$ \\
\\
\hline
4m & m+1 & m & m & m-1\\
4m+1 & m+1 & m & m & m\\
4m+2 & m+1 & m+1 & m & m\\
4m+3 & m+1 & m+1 & m+1 & m\\
\hline
\end{tabular}
\end{center}
%\label{}
\end{table}%

Later for the numerical construction we will need the next Proposition.
\begin{proposition}\label{P:250203-3}
Let $A$ be a subset of $\Z/M\Z$ and $f$ be a function in $E_M$ such that $f(a)=0$ for $a\in A$. Assume that $\widehat f(b)=f(b)$ for $b\not\in A$, then $\widehat f=f$, that is $\widehat f(a)=0$ for $a\in A$. 
\end{proposition}
\begin{proof}
Since $F$ is an isometry, we have
\[\sum_{a\in\Z/M\Z}|\widehat f(a)|^2=\sum_{a\in\Z/M\Z}|f(a)|^2=\sum_{b\notin A}|f(b)|^2=\sum_{b\not\in A}|\widehat f(b)|^2\]
Therefore 
\[\sum_{a\in A}|\widehat f(a)|^2=0.\]
and all these coordinates are $0$.
\end{proof}

When $M=N^2$ where $N$ is an odd natural number, we have $N^2=4K+1$, the space $\Z/\N^2\Z$ can be identified with $\{-2K, -2K+1,\dots, -1,0,1,\dots 2K\}$ and $E_{N^2}$ with functions on this set. Therefore, the values $f(n)$ for $|n|\le 2K$  determines $f$, and we may speak of symmetric $f$, when $f(n)=f(-n)$.

\begin{proposition}\label{P:250131-1}
Let $N$ be an odd natural number and $T>0$ a real number with $N^2>4NT+1$. Then there is a non null, real and symmetric $f\in E_{N^2}$, with $\widehat f=f$ and vanishing for $|n|\le NT$. 
\end{proposition}
\begin{proof}
For $\lambda^4=1$ let $V_\lambda\subset E_{N^2}$ be the eigenspace corresponding to $\lambda$. Let $W\subset E_{N^2}$ be the space of symmetric functions vanishing for $|n|\le NT$. For $f$ symmetric we have $f(n)=f(-n)=(F^2 f)(n)$ therefore $F^2f=f$, this implies that $W\subset V_1\oplus V_{-1}$. 

Each $w\in W$ can be written in a unique way as $w=u+v$ with $u\in V_1$ and $v\in V_{-1}$. Therefore there is a linear map $\pi\colon W\to V_{-1}$ that sends $w$ into $v$.

Since $N^2=4K+1$, we have  $\dim(V_{-1})=K$ and $\dim(W)=2K-\lfloor NT\rfloor$, so that $\dim(W)>\dim(V_{-1})$ is equivalent to $K>\lfloor NT\rfloor $ with $K=\frac{N^2-1}{4}$ and we assume that $N^2>4NT+1$.

It follows that $\pi$ cannot be injective and there are two points 
$w=u+v$ and $w'=u'+v$ with $w\ne w'$. Then $w-w'=u-u'\in W\cap V_1$. Therefore, the element $f=w-w'$ satisfies the conditions in our theorem, except that it is real. But for a symmetric function $f$ we have 
\begin{equation}\label{E:250203-1}
\widehat f(n)=\frac{1}{N}\sum_{m=-2K}^{2K} f(m)e^{-2\pi i \frac{mn}{N^2}}=\frac{f(0)}{N}+\frac{2}{N}\sum_{n=1}^{2K}f(n)\cos\Bigl(2\pi\frac{nm}{N^2}\Bigr).
\end{equation}
It follows that the real and imaginary part of $f$ are transformed independently one of the other. If $f\ne0$ is not real, one of the functions $\Re f$ or $\Im f\ne0$ is a real functions satisfying all our conditions.
\end{proof}

\begin{proposition}\label{P:250131-4}
Let $N$ be an odd natural number, and $T>0$ with $N^2>4NT+1$. There exist a non trivial crystalline measure 
\begin{equation}\label{E:250131-3}
\mu=\sum_{n\in\Z}c_n\delta_{n/N}
\end{equation}
that is real ($c_n\in\R$), symmetric ($c_{-n}=c_n$), periodic ($c_{n+N^2}=c_n$) with $\widehat\mu=\mu$  and such that the restriction of $\mu$ to $[-T,T]$ is $0$ ($c_n=0$ for $|n|\le NT$). 
\end{proposition}
\begin{proof}
Let $f$, the function constructed in Proposition \ref{P:250131-1}. Let us define $c_n=f(n)$ for $|n|\le(N^2-1)/2$ and extend $c_n$ periodically by $c_{n+N^2}=c_n$. 
Then by Proposition \ref{P-250131-2} the measure defined by \eqref{E:250131-3} satisfies $\widehat\mu=\mu$ because $\widehat f=f$. All the other properties of $\mu$ reduce to the properties of the values of $f(n)$. Since $f\ne0$, we have $\mu\ne0$.
\end{proof}

\section{Zeta function sharing many properties with Riemann zeta}\label{S:newzeta}

For $a>0$ let  $\zeta(s,a)$ denote the Hurwitz zeta function 
\[\zeta(s,a)=\sum_{n=0}^\infty\frac{1}{(n+a)^s}.\]
This functions extends to a meromorphic function on $\C$ with a unique simple pole at $s=1$ with residue $1$. (For the Hurwitz function see \cite{A}).

We need also the bound. There exist constants $C$ and $c$ depending on $a$, such that 
\begin{equation}\label{E:250201-5}
|(s-1)\zeta(s,a)|\le C e^{c|s|\log(2+|s|)}.
\end{equation}
%  Proof in my file \texttt{250324-Hurwitz.pdf}

\begin{proposition}\label{P:250201-1}
Let $N$ be an odd natural number and $T\ge1$ a natural number such that $N^2>4NT+1$. 
There exists an entire function $g_N(s)$ such that
\begin{itemize}
\item[(a)] For $\sigma>\sigma_0$ it is given by an absolutely convergent Dirichlet series \[g_N(s)=\sum_{n=1}^\infty a_n \lambda_n^{-s},\qquad T<\lambda_1<\lambda_2<\cdots\]
\item[(b)] $g_N(s)$ satisfies the functional equation of zeta
\begin{equation}\label{E:funceq}
\pi^{-s/2}\Gamma(\tfrac{s}{2})g_N(s)=\pi^{-\frac{1-s}{2}}\Gamma(\tfrac{1-s}{2})g_N(1-s).\end{equation}
\item[(c)] $N^sg_N(s)$ is a finite linear combination of Hurwitz zeta functions
\begin{equation}\label{E:250131-5}
g_N(s)=\frac{1}{N^s}\sum_{r=NT+1}^{(N^2-1)/2} c_r \zeta(s,r/N^2).
\end{equation}
\end{itemize}
\end{proposition}
\begin{proof}
Let $\mu$ the measure constructed in Proposition \ref{P:250131-4}. Since $\widehat\mu=\mu$ Theorem \ref{T:cryszeta} implies that the Dirichlet series 
\begin{equation}\label{E:250201-4}
g_N(s):=\sum_{n=1}^\infty \frac{c_n}{(n/N)^s},
\end{equation}
extends to a meromorphic function with a simple pole at $s=1$ with residue $c_0$, that satisfies the functional equation of zeta \eqref{E:funceq}. 

In our construction $c_n=0$ for $|n|\le NT$. Therefore, the residue $c_0=0$ and $g_N(s)$ is an entire function. The first non null term of the Dirichlet series is for $n=NT+1$, so it starts with $c_{NT+1}(T+1/N)^{-s}$, that is $\lambda_1=T+1/N>T$ as required. 

We have 
\[\sum_{n=1}^\infty \frac{c_n}{(n/N)^s}=\sum_{r=NT+1}^{(N^2-1)/2} c_r\sum_{m=0}^\infty \frac{1}{((mN^2+r)/N)^s}=\sum_{r=NT+1}^{(N^2-1)/2} \frac{c_r}{N^s}\sum_{m=0}^\infty\frac{1}{(m+r/N^2)^s}\]
since $c_r=0$ for $0\le r\le NT$, and the periodicity $c_{n+N^2}=c_n$. This is \eqref{E:250131-5}.
\end{proof}

\begin{proposition}\label{P:250201-2}
In the conditions of Proposition \ref{P:250201-1} let $\zeta_N(s)=\zeta(s)+\delta g_N(s)$ where $\delta$ is a real number. Then
\begin{itemize}
\item[(a)] $\zeta_N(s)$ is meromorphic with a unique simple pole at $s=1$ with residue $1$. 
\item[(b)] $\zeta_N(x)\in\R$ for $x$ real $x\ne1$.
\item[(c)] $\xi_N(s):=\frac{s(s-1)}{2}\Gamma(s/2)\pi^{-s/2}\zeta_N(s)$ satisfies the functional equation $\xi_N(s)=\xi_N(1-s)$ and it is real for $\sigma=\frac12$.
\end{itemize}
\end{proposition}
\begin{proof}
Since $g_N(s)$ is entire (a) follows from the properties of $\zeta(s)$. In equation \eqref{E:250131-5} the coefficients $c_n$ are real numbers, so  $g_N(s)$ is real for real $s$. The two functions $\zeta(s)$ and $g_N(s)$ satisfies the same functional equation. Therefore, $\xi_N(s)=\xi_N(1-s)$. By (b) the function $\xi_N(s)$ is real for real $s$, therefore $\xi_N(\overline s)=\overline{\xi_N(s)}$. For $s=\frac12+it$ we have $\xi_N(s)=\xi_N(1-s)=\overline{\xi_N(1-\overline s)}$. That is $\xi_N(s)$ take complex conjugate values at symmetric points on the critical line. This gives us that it takes real values on the critical line.
\end{proof}

We want to construct a function $\Xi_N(t)$ so that its zeros is a Riemann sequence.
Conditions (a) and (c) of the definition of Riemann sequence, are related to certain region without zeros for $\Xi_N(t)$. This is obtained in the next Proposition. 

\begin{proposition}\label{P:250202-1}
Assume that $N\ge2$, and let $\Xi_N(t):=\xi_N(\frac12+it)$, where $\xi_N(s)$ is defined in Proposition \ref{P:250201-2}.

There is a $\delta_0>0$ such that for $|\delta|<\delta_0$, the function $\Xi_N(t)$ do not vanish for $|\Im(t)|\ge3/2$, nor on the square $[-3/2,3/2]\times[-3/2,3/2]$. 
\end{proposition}
\begin{proof}
Since $\Xi(t)$ is even we have only to prove that $\Xi_N(t)\ne0$ in $\Im(t)\le -3/2$, nor in the rectangle $[-3/2,3/2]\times[-3/2,0]$. This is equivalent to say that $\xi_N(s)$ do not vanish for $\sigma\ge2$ or the rectangle $R:=[1/2,2]\times[-3/2,3/2]$.

By  \eqref{E:250201-4} we have for $\sigma\ge2$
\[|\zeta_N(s)|=\Bigl|\sum_{n=1}^\infty\frac{1}{n^s}+\delta\sum_{n=NT+1}^\infty \frac{c_n}{(n/T)^s}\Bigr|\ge 1-(\zeta(2)-1)-\delta\sum_{n=NT+1}^\infty\frac{|c_n|}{(n/T)^2}>0,\]
taking $|\delta|<\delta_0$ for an adequate $\delta_0$. Therefore $\xi_N(s)$ do not vanish for $\sigma\ge2$.

The function $(s-1)\zeta(s)$ is continuous and do not vanish on the compact set  $R$, so it has a minimum value $|(s-1)\zeta(s)|\ge m>0$ for $s\in R$. The function $(s-1)g_N(s)$ is continuous on this compact set so it has a maximum $|(s-1)g_N(s)|\le M$  for $s\in R$.
Then we have $|(s-1)\zeta_N(s)|\ge m-\delta M>0$ if $|\delta|\le m/2M$.
\end{proof}

\begin{proposition}
For $N\ge2$ and $|\delta|\le \delta_0$ the equality
\begin{equation}
\Xi_N(t)=\Xi_N(0)\prod_{\Re\alpha>0}\Bigl(1-\frac{t^2}{\alpha^2}\Bigr),
\end{equation}
holds true. Here $\alpha$ run through the zeros of $\Xi(t)$ with $\Re\alpha>0$.
\end{proposition}
\begin{proof}
Since the Hurwitz zeta function satisfies the bound \eqref{E:250201-5}, we have
$|(s-1)\zeta_N(s)|\le Ce^{c|s|\log(2+|s|)}$ for some positive constants $C$ and $c$.
The function $\xi_N(s)$ is the product of this function with $\frac{s}{2}\pi^{-s/2}\Gamma(s/2)$, For $\sigma\ge1/2$ we have $|\Gamma(s)|\le Ce^{|s|\log(2+|s|)}$. It follows that 
\[|\Xi(t)|\le C e^{c|t|\log(2+|t|)}.\]
By Jensen theorem (see \cite{B}*{1.2} this implies that the number $N(r)$ of zeros of $\Xi_N(t)$ in the disc of center $0$ and radius $r>0$ is bounded by $Cr\log (r+2)$.
Then $\sum|\alpha|^{-2}$ extended to the zeros $\alpha$ of $\Xi_n$ is bounded.
By Weierstrass theorem we will have
\[\Xi_N(t)=e^{at+b}\prod_\alpha\Bigl(1-\frac{t}{\alpha}\Bigr)e^{-t/\alpha}.\]
But $\Xi_N(t)=\Xi_N(-t)$ mean that for each root $\alpha$ of $\Xi_N(t)$  the function vanishes also on $-\alpha$. Joining the corresponding terms we get $\Xi_N(t)=e^{at+b}\prod_{\Re\alpha>0}(1-\frac{t^2}{\alpha^2})$. Again $\Xi_N(t)=\Xi_N(-t)$ implies $a=0$ and $e^b=\Xi(0)$.
\end{proof}

\begin{theorem} \label{T:mainconst}
Let $N$ be and odd natural number and $T\ge2$ a natural number, with $N^2>4NT+1$. Then  there is a $\delta_0>0$ such that for $|\delta|<\delta_0$ the sequence $(\alpha_n)$ of the zeros of $\Xi_N(s)$ with $\Re\alpha>0$, repeated with its multiplicity and ordered so that 
$0<\Re\alpha_1\le \Re\alpha_2\le \cdots$ is a Riemann sequence.
\end{theorem}

\begin{proof}
The function $\zeta_N(s)$ share with $\zeta(s)$ the following properties:
\begin{itemize}
\item[(1)] Is meromorphic in $\C$ with a unique simple pole at $s=1$ with residue $1$.
\item[(2)] $\zeta_N(s)\simeq 1$.
\item[(3)] It has the product representation
\[\zeta_N(s)=\Xi_N(0)\frac{\pi^{s/2}}{\frac12s(s-1)\Gamma(\frac{s}{2})}\cdot\prod_{n_1}^\infty \Bigl(1+\frac{(s-\frac12)^2}{\alpha_n^2}\Bigr).\]
\end{itemize}
With these properties and reasoning as in Section \ref{suno} we show that the sum $\sum\frac{2z}{z^2+\alpha^2}$ has the asymptotic expansion \eqref{mainalpha}. This is property (e) of the definition of Riemann sequence.

The terms of the sequence are the zeros $\alpha$ of $\Xi_N(t)$ with $\Re(\alpha)>0$. Proposition \ref{P:250202-1} shows that $|\Im(\alpha)|<3/2$ (this is condition (c) in Definition of Riemann sequence) and that $\Re(\alpha)>3/2$ (this is stronger than condition (a)). Between the two shows that $|\Im(\alpha)|<\Re(\alpha)$ (and this is condition (d) of the Definition). Condition (b) follows from the fact that $\Xi_N$ is real on the real axis.
\end{proof}

\begin{remark}
There is no reason for the function $\zeta_N(s)$ to have an Euler product. So we do not expect the zeros of this function to be on the critical line. In the next section we construct a concrete example and show numerically that it does not satisfies the Riemann hypothesis. This proves that there are \textbf{complex} Riemann sequences. 
Our main remaining question is whether there is \textbf{real} Riemann sequences. Given that the zeros of zeta off the critical line (in case they exist) are very far. To show that there is not real Riemann sequence, may be the best option to disprove the Riemann hypothesis.
\end{remark}

\begin{remark}
In case there is some real Riemann sequence it will be interesting to prove the uniqueness. We may use some other functional equations satisfied by several Dirichlet series with Euler products. In this case there is no uniqueness, but we can expect uniqueness in the case of the Riemann zeta function. 
\end{remark}

\section{Numerical example of Complex Riemann Sequence}\label{S:numeric}

In the previous section we have shown the existence of Riemann sequences. Since they are defined from Dirichlet series that do not have an Euler product it is likely that their zeros do not satisfy the Riemann hypothesis, i.e. the sequences obtained are not real. We now construct a concrete example of a complex Riemann sequence. Thus demonstrating that there exist non-real Riemann sequences.

We can construct examples with additional properties. For example, in that the Dirichlet series starts with $1+2^{-s}$ such that the lines where $\zeta_\alpha(s)$ are real or pure imaginary are in bijective correspondence with the corresponding lines for $\zeta(s)$. Or so that the residue at $s=1$ of our function $\zeta_\alpha(s)$ is equal to $1$. But we have preferred to give an example that is easy to define and does not depend on complicated constants. 

Define the measure  $\sigma=\frac13\sum_{n\in\Z}\delta_{n/3}$.  To simplify the notation let $a=\sqrt{3}$ and define the trigonometric polynomial
\[
P(x)=(1-a)\cos(2\pi x)+2\cos(5\pi x/2)-(1+a)\cos(3\pi x)+2\cos(7\pi x/2)+(1-a)\cos(4\pi x).
\]
Then define the measures
\[\mu=P\sigma,\quad \nu=\widehat{\mu}=\widehat P*\widehat\sigma=\widehat P * \sum_{n\in\Z}\delta_{3n}.\]
An easy computation yields then 
\[\nu=(1-a)\sum_{n\in\Z}(\delta_{3n+1}+\delta_{3n+2})+2\sum_{n\in\Z}(\delta_{3n+5/4}+\delta_{3n+7/4})-(1+a)\sum_{n\in\Z}\delta_{3n+3/2}.\]
We need also a more explicit expression for $\mu$
\[\mu=P\sigma=\frac{P}{3}\sum_{n\in\Z}\delta_{n/3}=\sum_{n\in\Z}\frac{P(n/3)}{3}\delta_{n/3}\]
The expression for $P$ implies that $P(x+4)=P(x)$, it follows that putting $n=12k+r$ for $0\le r\le11$ we get the expression
\[\mu=\sum_{r=0}^{11}\frac{P(r/3)}{3}\sum_{k\in\Z}\delta_{4k+r/3}.\]
The values of the coefficients $c_r:=P(r/3)/3$ are as follows (with $a=\sqrt{3}$)
\[c_0=\frac{5-3a}{3},\quad c_1=c_2=c_{10}=c_{11}=0,\quad c_3=c_9=\frac{3-a}{3},\quad c_4=c_8=-\frac43,\]\[c_5=c_7=\frac{4a}{3},\quad c_6=-(1+a)\]

We are now in position to define our crystalline measure. The measure $\mu+\nu$ is real, symmetric, its proper Fourier transform, and its restriction to $[-1,1]$ is equal to 
\[\frac{3-a}{3}\delta_{-1}+(1-a)\delta_{-1}+\frac{5-3a}{3}\delta_0+\frac{3-a}{3}\delta_1+(1-a)\delta_1.\]
Since we have the coefficient of $\delta_1$ to be $1$, we multiply this measure by $(2-4/a)^{-1}$. 

The associated zeta function may be expressed in terms of the Hurwitz zeta function 
\begin{align*}
\zeta_M(s)&=-\frac{1+a}{2}\frac{1}{4^s}(\zeta(s,1/4)+\zeta(s,3/4))+\frac{6+4a}{3}\frac{1}{4^s}(\zeta(s,1/3)+\zeta(s,2/3))\\&-2(2+a)\frac{1}{4^s}(\zeta(s,5/12)+\zeta(s,7/12))+\frac{9+5a}{2}\frac{1}{4^s}\zeta(s,1/2)+\frac{3-a}{6}\frac{1}{4^s}\zeta(s,1)\\&+ \frac{3+a}{2}\frac{1}{3^s}(\zeta(s,1/3)+\zeta(s,2/3))-(3+2a)\frac{1}{3^s}(\zeta(s,5/12)+\zeta(7/12))+\frac{9+5a}{2}\frac{1}{3^s}\zeta(s,1/2).
\end{align*}
%This can be further simplified into 
%\begin{equation}
%\begin{aligned}
%\zeta_M(s)&=\Bigl(1+\frac{6+4a}{3}\frac{1}{(4/3)^s}+\frac{9+5a}{2}\frac{1}{(3/2)^s}+\frac{5+3a}{2^s}-\frac{6+3b}{3^s}-\frac{6+4a}{4^s}\Bigr)\zeta(s)\\
%&-\Bigl(\frac{3+2a}{3^s}+\frac{4+2a}{4^s}\Bigr)\bigl(\zeta(s,5/12)+\zeta(s,7/12)\bigr)).
%\end{aligned}
%\end{equation}

The first term of the Dirichlet series are
\[\zeta_M(s)=1-\frac{3+2a}{(5/4)^s}+\frac{2+\frac{4}{a}}{(4/3)^s}+\frac{\frac{9+5a}{2}}{(3/2)^s}-\frac{2(2+a)}{(5/3)^s}-\frac{3+2a}{(7/4)^s}+\frac{3(2+a)}{2^s}+\cdots\]

Since $\widehat{\mu+\nu}=\mu+\nu$ and $\mu+\nu$ is symmetric, the function $\zeta_M(s)$ extends to a meromorphic function on the complex plane, satisfies the functional equation $\pi^{-s/2}\Gamma(s/2)\zeta_M(s)=\pi^{-(1-s)/2}\Gamma((1-s)/2)\zeta_M(1-s)$ and has a unique simple pole at $s=1$ with residue $\frac{3-a}{6}$. Since $\mu+\nu$ is real the zeros of $\zeta_M(s)$ has the same symmetries as the ones of $\zeta(s)$. To show that the zeros of $\zeta_M(s)$ leads to a Riemann sequence, we need to show that the first complex zeros has ordinates greater than 1 and that they are contained in the strip $1-\sigma_0<\sigma\le \sigma_0$ for some $\sigma_0>1$. 

We want to determine $\sigma_0$ such that $\zeta_M(s)\ne0$ for $\sigma\ge\sigma_0$. 
The function $\zeta_M(s)$ is a Dirichlet series. The constant term appears only on
\[ -\frac{1+a}{2\cdot4^s}\zeta(s,1/4)+\frac{3+a}{2\cdot3^s}\zeta(s,1/3)=1+\cdots\]
Then the $\sigma_0$ is the point where 
\[1=\frac{1+a}{2}\Bigl(\frac{\zeta(\sigma,1/4)}{4^s}-1\Bigr)+\frac{3+a}{2}\Bigl(\frac{\zeta(\sigma,1/3)}{3^s}-1\Bigr)+\frac{1+a}{2\cdot4^s}\zeta(\sigma,3/4)+\cdots\]
where the dots represents all the other terms in $\zeta_M(s)$ taken with coefficients in absolute value

This happens for $\sigma_0=10.564029176912431172$. 
\begin{figure}[H]
\begin{center}
\includegraphics[height=0.95\vsize]{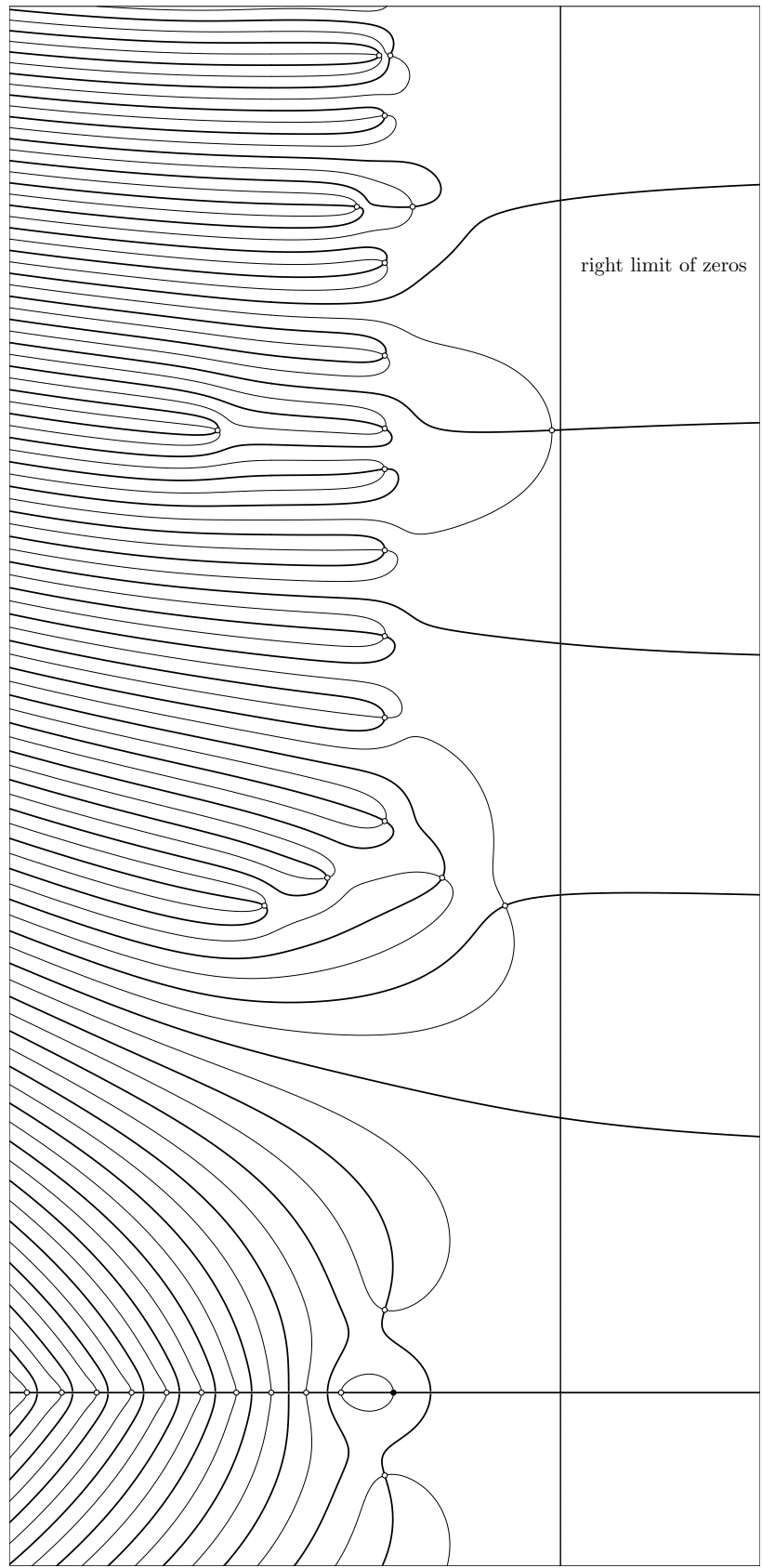}
\caption{x-ray of $\zeta_M(s)$ in $R=(-21,22)\times(-10,80)$.}
%\label{default}
\end{center}
\end{figure}

To show that conditions (a) and (d) of the definition of a Riemann sequence are satisfied we show an x-ray of the function $\zeta_M(s)$. This means we draw on a rectangle $R$ the lines where the function $\zeta_M(s)$ take real values (thick lines), and the lines where $\zeta_M(s)$ take purely imaginary values (thin lines). These lines cut at zeros and poles of $\zeta_M(s)$.

We may understand the x-ray as computing the variation of the argument along the boundary of the rectangle. We know in this way the number of zeros of $\zeta_M(s)$ contained in this rectangle. In our case 31 zeros and a pole. We then compute the zeros, so that we know the computed zeros are all contained in this rectangle. The complex zeros with positive imaginary part are

\[\begin{array}{ll|ll}
\beta & \gamma & \beta & \gamma\\
0.5 & \phantom{2}4.7753735547 & 0.5 & 53.2934095839 \\
-6.3939983623& 28.0995236414 & \,10.0731303207 & 55.5328071355 \\
\phantom{-}7.3939983623 & 28.0995236414 & -9.0731303207 & 55.5328071355\\
-2.7891403465 & 29.7107114647 & 0.5 & 55.6380182916 \\
\phantom{-}3.7891403465 & 29.7107114647 & 0.5 & 59.8440884874 \\
0.5 & 32.9826068738 & 0.5 & 65.1982600562\\
0.5 & 38.9509449796 & 0.5 & 73.6917293006\\
0.5 & 43.6565105315 & 0.8205883718 & 77.1648164218\\
0.5 & 48.6090990060 & 0.1794116281 &  77.1648164218
\end{array}\]

\section{Properties of Riemann sequences}\label{sproperties}  

In what follows we assume that we have a Riemann sequence $(\alpha)$  and try to get information about it. 

\begin{proposition}\label{simpleboundN}
Let $(\alpha_n)$ be a Riemann sequence, and for $x>0$ let $N(x)$ be the number of terms with $\Re\alpha_n\le x$, then $N(x)\le Cx\log (x+1)$ for some  constant $C$, and  therefore $\sum_\alpha |\alpha|^{-1-\varepsilon}<+\infty$ for any $\varepsilon>0$.
\end{proposition}
\begin{proof}
For $x> 1$ and $\alpha=\gamma+i\beta$ with $|\beta|<C$ (by condition (c) of the definition of Riemann sequences) and $\gamma\le x$. When $\beta=0$
\[\frac{2x}{x^2+\gamma^2}\ge \frac{1}{x},\]
and when $\beta>0$, with $\beta\le\gamma$ by condition (d) of the definition of Riemann sequence 
\[\frac{2x}{x^2+\alpha^2}+\frac{2x}{x^2+\overline{\alpha}^2}=2\frac{2x(x^2+\gamma^2-\beta^2)}{(x^2+\gamma^2-\beta^2)^2+4\beta^2\gamma^2}\ge2\frac{2x}{5(x^2+\gamma^2-\beta^2)}\ge 2\frac{2x}{5x^2}\ge 2\frac{2}{5x}.\]
Therefore for $x>x_0$ we have
\[N(x)\frac{2}{5x}\le \sum_{\alpha}\frac{2x}{x^2+\alpha^2}\le \log x.\]
Hence for $x>x_0$ we have $N(x)\le 3x\log x$, and then there is a constant with $N(x)\le Cx\log (x+1)$ for $x>1$.
\end{proof}

It follows that the series in the definition of a Riemann sequence converges uniformly on compact sets of  $\C\smallsetminus\{\pm i \alpha_n\colon n\in \N\}$, and defines a meromorphic function with poles at the points $\pm i\alpha_n$.

\begin{theorem} Let $(\alpha_n)$ be a Riemann sequence. 
There exists a real  constant $B$ such that
\begin{equation}\label{eq17}
\sum_{\alpha}\log\Bigl(1+\frac{z^2}{\alpha^2}\Bigr)\simeq
\frac{z}{2}\log\frac{z}{2\pi}-\frac{z}{2}+\frac{7}{4}\log z+B-
\sum_{n=1}^{\infty}\frac{a_{n+1}}{n z^n},
\end{equation}
\end{theorem}

\begin{proof}
It is well known that from \eqref{mainalpha} it follows
\begin{equation}\label{E:28}
-\int_z^\infty\Bigl\{\sum_{\alpha}\frac{2t}{t^2+\alpha^2}-\frac12\log\frac{t}{2\pi}
+\frac{a_1}{t}\Bigr\}\,dt\simeq -\sum_{n=2}^\infty \frac{a_n}{(n-1)z^{n-1}}
\end{equation}
But the left hand side of \eqref{E:28} and
\begin{equation}
\sum_{\alpha}\log\Bigl(1+\frac{z^2}{\alpha^2}\Bigr)-
\frac{z}{2}\log\frac{z}{2\pi}+\frac{z}{2}+a_1\log z
\end{equation}
are two primitives of the integrand, so that they differ on  a constant $B$. Therefore
\begin{equation*}
\sum_{\alpha}\log\Bigl(1+\frac{z^2}{\alpha^2}\Bigr)-
\frac{z}{2}\log\frac{z}{2\pi}+\frac{z}{2}+a_1\log z-B\simeq-
\sum_{n=1}^{\infty}\frac{a_{n+1}}{n z^n}.
\end{equation*}

From the above we get 
\begin{equation}
B=\lim_{x\to+\infty}\sum_{\alpha}\log\Bigl(1+\frac{x^2}{\alpha^2}\Bigr)-
\frac{x}{2}\log\frac{x}{2\pi}+\frac{x}{2}-\frac74\log x.
\end{equation}
This combined with the fact that the multiplicities satisfies $\nu(\alpha)=\nu(\overline{\alpha})$ implies that $B$ is real.
\end{proof}

Given a Riemann sequence $(\alpha_n)$ we are going to define two associated functions $\zeta_\alpha(s)$ and $\Xi_\alpha(t)$ similar to the 
functions $\zeta(s)$ and $\Xi(s)$. We define $\Xi_\alpha(0)$ by the equation 
\begin{equation}
B=\frac14\log\frac{\pi}{2}-\log\Xi_\alpha(0)
\end{equation} 
(compare the value of $A$ in Theorem 
\ref{one}). Then define
\begin{align}\label{defxizeta}
\Xi_\alpha(t)&=\Xi_\alpha(0)\prod_{\alpha}\Bigl(1-\frac{t^2}{\alpha^2}\Bigr),\\
\zeta_\alpha(s)&=\Xi_\alpha(0)\frac{\pi^{\frac{s}{2}}}{(s-1)\Gamma(1+s/2)}\prod_\alpha
\Bigl(1+\frac{(s-1/2)^2}{\alpha^2}\Bigr)\label{defzeta}
\end{align}

These functions depend on the sequence $(\alpha)$ and when $(\alpha)=(\tau)$ we get the familiar $\zeta(s)$ and $\Xi(t)$.

\begin{theorem}
$\Xi_\alpha(t)$ is an even entire function 
of order  $1$.  Its zeros are the numbers $\alpha_n$ and $-\alpha_n$. 
$\Xi_\alpha(t)$ is real for $t$ real and for $t$ purely imaginary. 
\end{theorem}

\begin{proof}
Since $(\alpha_n)$ is a Riemann sequence $\sum_\alpha|\alpha|^{-2}<\infty$. Hence the product in \eqref
{defxizeta}  converges uniformly on compact sets and define an entire function.

Let $|t|\le r$,  is clear that for each factor with $\alpha\in\R$, we have $\Bigl|1-\frac{t^2}{\alpha^2}\Bigr|\le (1+\frac{r^2}{\alpha^2})$ when $\alpha$ is real. In the case of a complex $\alpha$, we may 
associate the factors for $\alpha$ and $\overline{\alpha}$. By condition (d) of the definition of Riemann sequence $|\arg(\alpha_n)|<\pi/4$ and therefore $\Re(\alpha_n^{-2})>0$. It follows that for $n\ge n_0$, 
\begin{multline*}
\Bigl|\Bigl(1-\frac{t^2}{\alpha^2}\Bigr)\Bigl(1-\frac{t^2}{\overline{\alpha}^2}\Bigr)
\Bigr|=\Bigl|1+\frac{t^4}{|\alpha|^4}-2t^2\Re\frac{1}{\alpha^2}\Bigr|\le 
1+\frac{r^4}{|\alpha|^4}+2r^2\Re\frac{1}{\alpha^2}=\\
=\Bigl(1+\frac{r^2}{\alpha^2}\Bigr)\Bigl(1+\frac{r^2}{\overline{\alpha}^2}\Bigr).
\end{multline*}

Then for $|t|\le r$
\[|\Xi(\alpha)|=\Bigl|\Xi_\alpha(0)\prod_\alpha\Bigl(1-\frac{t^2}{\alpha^2}\Bigr)\Bigr|\le |\Xi_\alpha(0)|\prod_\alpha\Bigl(1+\frac{r^2}{\alpha^2}\Bigr)=|\Xi(ir)|.\]

It follows that $|\Xi_\alpha(t)|\le  |\Xi_\alpha(ir)|$. Then by \eqref{eq17} we get 
\begin{equation}\label{boundlogM}
\log M(r)\sim \frac12 r\log r\qquad r\to+\infty.
\end{equation}
This proves that the order of $\Xi_\alpha(t)$ is $1$. 

By  definition the function $\Xi_\alpha(t)$ is even.  The 
order of multiplicity  $\nu(\alpha)=\nu(\overline{\alpha})$, hence the factors in the definition of $\Xi_\alpha(t)$ may be paired so that the product is real for $t$ real or 
purely imaginary, since we know that $\Xi_\alpha(0)$ is real we get the reality of 
$\Xi(x)$ and $\Xi(ix)$ for $x\in\R$.
\end{proof}

\begin{corollary}\label{expconv}
The exponent of convergence of a Riemann sequence $(\alpha)$ is $1$. The series 
$\sum_n|\alpha_n|^r$ converges for all $r>1$ but diverges for $0\le r\le1$.
\end{corollary}

\begin{proof} Since the order of $\Xi_\alpha(t)$ is $1$ the series converges for all $r>1$. Being $\Xi_\alpha(t)$ an even function with  $\int_2^\infty r^{-2}\log M(r)\,dr=+\infty$,  it follows that $\sum_n|\alpha_n|^{-1}=\infty$. (Cfr. Boas \cite{B}*{Thm.~2.11.1 applied to $\Xi_\alpha(\sqrt{t})$}.)
\end{proof}

\begin{theorem}
$\zeta_\alpha(s)$ is a meromorphic function on $\C$ with a unique simple pole at the 
point $s=1$ with residue 
\begin{equation}
\Res_{s=1}\zeta_\alpha(s)=2\,\Xi_\alpha(0)\prod_\alpha\Bigl(1+\frac{1}{4\alpha^2}\Bigr).
\end{equation}
Its zeros are the even negative integers $-2$, $-4$, $-6$, \dots and the numbers
$\tfrac12+i\alpha$, $\tfrac12-i\alpha$. These last zeros (the non trivial zeros) are 
contained on the `critical strip' $|\sigma-\frac12|<C$.

Also $\zeta_\alpha(s)$ satisfies the same  functional equation as the Riemann zeta function. 
\begin{equation}
\pi^{-\frac{s}{2}}\Gamma(s/2)\zeta_\alpha(s)=\pi^{-\frac{1-s}{2}}\Gamma
\bigl((1-s)/2\bigr)\zeta_\alpha(1-s).
\end{equation}
\end{theorem}

\begin{proof}
All the assertions are easily obtained from the definition \eqref{defxizeta}. The functional equation is equivalent to the assertion that $\Xi_\alpha(t)$ is an even function. 
\end{proof}

\begin{theorem}
We have
\begin{equation}
\log\zeta_\alpha(s)\simeq 0,\qquad \zeta_\alpha(s)\simeq1.
\end{equation}
\end{theorem}

\begin{proof}
Taking log on the definition \eqref{defzeta} and applying \eqref{defxizeta} we get 
\begin{multline}\label{logonalldos}
\log\zeta_\alpha(\tfrac12+z)=\\
\sum_\alpha\log\Bigl(1+\frac{z^2}{\alpha^2}\Bigr)+
\Bigl(\frac{z}{2}+\frac14\Bigr)
\log\pi-\log\Bigl(z-\frac12\Bigr)-\log\Gamma\Bigl(\frac{z}{2}+\frac{5}{4}\Bigr)+
\log\Xi_\alpha(0)
\end{multline}
this is similar to \eqref{logonall}. Now we may apply \eqref{eq7}, \eqref{eq8} and 
\eqref{eq17} to get exactly that $\log\zeta_\alpha(\tfrac12+z)\simeq0$. 

Observe that $\log\zeta_\alpha(\tfrac12+z)\simeq0$ is equivalent to say that on each angle
$|\arg z|\le \theta<\frac{\pi}{2}$ and for each natural number $N$ we have
\begin{equation}
\lim_{z\to\infty}z^N \log\zeta_\alpha(\tfrac12+z)=0
\end{equation}
It is clear that this implies that for each $N$ there is a value of $R>1$  such that 
\begin{equation}
|\arg z|\le \theta \quad\text{and}\quad |z|>R \qquad \Longrightarrow\qquad\bigl|z^N \log\zeta_\alpha(\tfrac12+z)|\le 1
\end{equation}
Then we will have $|\zeta_\alpha(\tfrac12+z)-1|\le C|z|^{-N}$ for $|\arg z|\le \theta$ and $|z|>R$.
It follows that $\zeta_\alpha(\tfrac12+z)\simeq1$

Finally $\zeta_\alpha(s)\simeq1$ and $\zeta_\alpha(\tfrac12+s)\simeq1$ are equivalent assertions.
\end{proof}

\section{The function \texorpdfstring{$Z_\alpha$.}{Za.}}

\begin{definition}
Given a Riemann sequence $(\alpha_n)$ and for $\Re z>0$ we define
\begin{equation}
f(z):=\sum_{n=1}^\infty e^{-\alpha_n^2z}.
\end{equation}
The function $f$ is holomorphic in $\Re z>0$.
\end{definition}

\begin{proof}
We must show that the series converges uniformly in compact sets of $\Re z>0$. 
Let $K$ such  a compact. There exists $\delta>0$ and $M<\infty$ 
such that $\Re z\ge \delta>0$ and $|z|\le M$ for $z\in K$. 

Put $\alpha_n=a_n+ib_n$,  its real and imaginary parts and $z=x+iy\in K$. Then we have
\begin{displaymath}
|e^{-\alpha_n^2z}|= e^{-(a_n^2-b_n^2)x+2a_nb_n y }\le e^{-a_n^2 \delta+Mb_n^2+2Ma_n|b_n|}
\end{displaymath}
By condition (c) of the definition of Riemann sequence $|b_n|\le C$ we get for $n$  big enough
\begin{displaymath}
|e^{-\alpha_n^2z}|\le c e^{-a_n^2 \delta+ C Ma_n}\le c' e^{-\frac{\delta}{2}a_n^2}.
\end{displaymath}
Finally the numerical series $\sum e^{-\frac{\delta}{2}a_n^2}<+\infty$. This follows from Proposition \ref{simpleboundN} by partial summation.
\end{proof}

\begin{proposition}
For $|\arg(t)|<\frac{\pi}{4}$ we have
\begin{equation}
\int_0^\infty f(x)e^{-xt^2}\,dx=\sum_{n=1}^\infty\frac{1}{t^2+\alpha_n^2}
\end{equation}
\end{proposition}

\begin{proof}
Condition (d) in Definition \ref{D:RS} implies  $\Re(\alpha_n^2)>0$, so that each integral \[\int_0^\infty e^{-\alpha_n^2x}e^{-xt^2}\,dx=(t^2+\alpha_n^2)^{-1}\] is convergent. 

Since $\Re(\alpha_n)\to+\infty$ and $|\Im(\alpha_n)|\le C$, there is some $n_0$ such that $\Re(\alpha_n)>2|\Im(\alpha_n)|$ for $n\ge n_0$. With the same notation as before $\alpha_n=a_n+ib_n$ and for $t=u+iv$ with $\delta=u^2-v^2>0$, we have 
\[|e^{-\alpha_n^2 x}e^{-xt^2}|=e^{-(a_n^2-b_n^2)x-(u^2-v^2)x}\le e^{-(\frac34a_n^2+\delta)x},\qquad n\ge n_0,\quad x>0.\]
Therefore,
\[\int_0^\infty\Bigl|\sum_{n\ge n_0}e^{-\alpha_n^2} e^{-xt^2}\Bigr|\,dx\le \sum_{n=n_0}^\infty\int_0^\infty e^{-(\frac34a_n^2+\delta)x}\,dx=\sum_{n\ge n_0}\frac{1}{\frac34 a_n^2+\delta}<+\infty.\]

So the integration term by term is allowed and we get for $\Re t^2>0$,
\begin{equation}
\int_0^\infty f(x)e^{-xt^2}\,dx=\sum_{n=1}^\infty \int_0^\infty e^{-\alpha_n^2x}e^{-xt^2}\,dx=\sum_{n=1}^\infty \frac{1}{t^2+\alpha_n^2}.
\end{equation}
The restriction $\Re t^2>0$ is equivalent to $|\arg(t)|<\frac{\pi}{4}$.
\end{proof}

Now we may get the behavior of $f(x)$ when $x>0$ tends to $0$. 

\begin{theorem}
When $x$ real and $x\to0^+$ we have the asymptotic expansion
\begin{equation}\label{expasym}
2\sum_\alpha e^{-\alpha^2x}\sim \frac{1}{4\sqrt{\pi x}}\log\frac{e^{-C_0}}
{16\pi^2 x}
+\sum_{n=0}^\infty \frac{a_{n+1}}{\Gamma(1+n/2)}x^{n/2},\qquad x>0,\quad x\to0.
\end{equation}
\end{theorem}

\begin{proof}
Combining two entries of the Table of Integral Transform of Bateman \cite{Ba}*{4.3(1) p.~137 and 4.6(9) p.~148} we obtain\footnote{Beware that the second formula in the Bateman tables is not correct. Gradshteyn and Ryzhik tables \cite{GR}*{17.13.92} gives the same erroneous value.}
the following integral
\begin{equation}
\int_0^{+\infty} \frac{1}{4\sqrt{\pi x}}\log\frac{e^{-C_0}}{16\pi^2 x}e^{-x t^2}
\,dx=\frac{1}{2t}\log\frac{t}{2\pi},\qquad |\arg(t)|<\frac{\pi}{4}.
\end{equation}

Hence for $N\ge1$ we have
\begin{multline}
\int_0^{+\infty}\Bigl\{\frac{1}{4\sqrt{\pi x}}\log\frac{e^{-C_0}}{16\pi^2 x}+
\frac74
+\sum_{n=1}^{N-1} \frac{a_{n+1}}{\Gamma(1+n/2)}x^{n/2}\Bigr\}e^{-xt^2}\,dx=\\
=\frac{1}{2t}\log\frac{t}{2\pi}+\frac{7}{4t^2}+\sum_{n=1}^{N-1} \frac{a_{n+1}}{t^{n+2}}.
\end{multline}
Put
\begin{equation} \label{defgN}
g_N(x):=2f(x)-\Bigl\{\frac{1}{4\sqrt{\pi x}}\log\frac{e^{-C_0}}{16\pi^2 x}+
\frac74
+\sum_{n=1}^{N-1} \frac{a_{n+1}}{\Gamma(1+n/2)}x^{n/2}\Bigr\}
\end{equation}
Then 
\begin{equation}
G_N(t):=\int_0^\infty g_N(x)e^{-xt^2}\,dx
\end{equation}
is defined for $|\arg(t)|<\frac{\pi}{4}$ and is equal to 
\begin{equation}
G_N(t)=\sum_{n=1}^\infty \frac{2}{t^2+\alpha_n^2}-\frac{1}{2t}\log\frac{t}{2\pi}-\frac{7}{4t^2}-\sum_{n=1}^{N-1} \frac{a_{n+1}}{t^{n+2}}.
\end{equation}
Therefore $G_N(t)$ extends to an analytic function on $\Re t>\frac12$ and verifies
\begin{equation}
G_N(t)\simeq -\sum_{N+1}^\infty\frac{a_n}{t^{n+1}}
\end{equation}

In particular for each $0<\theta<\frac{\pi}{2}$  
\[
G_N(t)\quad\text{extends to $\Re t>\tfrac12$ and } G_N(t)\sim -\frac{a_{N+1}}{t^{N+2}}
\quad\text{when $t\to\infty$ in $|\arg(t)|\le \theta$.}
\]
Changing variable $s:=t^2$, it follows that
for any $0<\varphi<\pi$ 
\begin{equation}\label{hyp}
\int_0^\infty g_N(x)e^{-sx}\,dx
\end{equation}
extends to an analytic function $H_N(s)$ on the region $V=\{(x,y):x+y^2>1/4\}$ and for any $0<\varphi<\pi$
\begin{equation}
H_N(s)\sim -\frac{a_{N+1}}{s^{1+\frac{N}{2}}} \quad\text{when $s\to\infty$ in $|\arg(s)|\le \varphi,\quad s\in V$.}
\end{equation}

Now we apply a classical Theorem  (see Doetsch \cite{Do}*{p.~503--504 Satz 3}) to get from \eqref{hyp} that 
\begin{equation}\label{ordergN}
g_N(x)\sim-\frac{a_{N+1}}{\Gamma(1+\frac{N}{2})}x^{\frac{N}{2}},\qquad x>0,\quad x\to0.
\end{equation}
But that this is true for all $N$ is just  the meaning of \eqref{expasym}.
\end{proof}

For $\Re s>1$ we define
\begin{equation}
Z_\alpha(s):=\sum_{n=1}^\infty \frac{1}{\alpha_n^s},\qquad \Re s>1.
\end{equation}

By definition of Riemann sequence $\Re\alpha_n\ge\Re\alpha_1>1$ and $|\arg(\alpha_n)|\le\frac{\pi}{4}$. Hence $\alpha_n^{-s}$ has a well defined sense, taking the principal value of $\log\alpha_n$.  It follows that  $|\alpha_n^{-s}|=e^{-\sigma\log|\alpha_n|+t\arg(\alpha_n)}\le e^{\frac{\pi}{4}|t|}|\alpha_n|^{-\sigma}$. So 
by Corollary \ref{expconv} the series $\sum_n\alpha_n^{-s}$ is uniformly convergent 
on compact sets of $\Re s>1$. 

It follows that $Z_\alpha(s)$ is a holomorphic function on $\Re s>1$. 

\begin{theorem}\label{T:intZ}
For $\Re s>1$ we have
\begin{equation}\label{modexp}
\Gamma(s/2)Z_\alpha(s)=\int_0^{+\infty}f(x)x^{\frac{s}{2}}\frac{dx}{x}.
\end{equation}
\end{theorem}

\begin{proof}
By definition of Riemann sequence $|\arg(\alpha_n)|<\pi/4$, therefore $\Re(\alpha_n^2)>0$ and each integral
\begin{equation}
\int_0^\infty e^{-\alpha_n^2 x}x^{\frac{s}{2}}\frac{dx}{x}=\frac{\Gamma(s/2)}{\alpha_n^s},\qquad \sigma>0,
\end{equation}
is convergent.
By condition (c) of the definition of Riemann sequence $|\Im(\alpha_n)|\le C$, then $\Re(\alpha_n)\to+\infty$ and we have $|\Im(\alpha_n)|\le \frac12\Re(\alpha_n)$ for $n\ge n_0$. It follows that $|\alpha_n|^2\le \frac54\Re(\alpha_n)^2$ and 
\[\Re(\alpha_n^2)=\Re(\alpha_n)^2-\Im(\alpha_n)^2\ge \tfrac34\Re(\alpha_n)^2\ge \frac35|\alpha_n|^2.\]
Therefore, for $\sigma>1$
\[\Bigl|\sum_{n\ge n_0}\int_0^\infty e^{-\alpha_n^2 x}x^{\frac{s}{2}}\frac{dx}{x}\Bigr|\le \sum_{n\ge n_0}\int_0^\infty e^{-\frac35|\alpha_n|^2 x}x^{\frac{\sigma}{2}}\frac{dx}{x}=\sum_{n\ge n_0}\frac{\Gamma(\sigma/2)}{(3/5)^{\sigma/2}|\alpha_n|^\sigma}<+\infty.\]
This justify the interchange of sum and integral and we get 
\begin{equation}
\int_0^\infty f(x)x^{\frac{s}{2}}\frac{dx}{x}=\sum_{n=1}^\infty
\int_0^\infty e^{-\alpha_n^2 x}x^{\frac{s}{2}}\frac{dx}{x}=\Gamma(s/2)Z_\alpha(s).\qedhere
\end{equation}
\end{proof}

\begin{theorem}
The function $Z_\alpha(s)$ extends to a meromorphic function on $\C$. 

The function $Z_\alpha(s)$ has a double pole at $s=1$ with residue $-\frac{\log2\pi}{2\pi}$ and main part
\begin{equation}
\frac{1}{2\pi(s-1)^2}-\frac{\log2\pi}{2\pi(s-1)}.
\end{equation}
It has simple poles at the odd negative integers with main part
\begin{equation}
(-1)^{n}\frac{  a_{2n}}{\pi(s+2n-1)} = (-1)^n\frac{(1-2^{1-2n})}{2\pi}\frac{B_{2n}}
{2n}\frac{1}{s+2n-1}.
\end{equation}
The function $Z_\alpha(s)$ has no other poles.

Finally we have 
\begin{equation}\label{Zvalue}
Z_\alpha(-2n)=(-1)^{n}\frac{a_{2n+1}}{2}=(-1)^n\frac{8-E_{2n}}{2^{2n+3}}.
\end{equation}
\end{theorem}

\begin{proof}
First notice that there is a constant $\delta>0$ such that $|f(x)|\le Ce^{-\delta x}$ for all $x>1$. As in the proof of Theorem \ref{T:intZ} there is some $n_0$ such that for $n\ge n_0$ $|e^{-\alpha_n^2x}|\le e^{-\frac35|\alpha_n|^2x}$, and for each $n<n_0$ we have $|e^{-\alpha_n^2x}|\le e^{-(\Re(\alpha_n)^2-\Im(\alpha_n)^2)x}=e^{-\delta_nx}$ with $\delta_n>0$. Since $\Re(\alpha_n)^2\to+\infty$ there is some $\delta>0$ such that it is the minimum of all these exponents. Therefore, for $x>1$
\begin{equation}\label{boundf}
|f(x)|\le e^{-\delta x}\Bigl(\sum_{n<n_0}e^{-(\delta_n-\delta)x}+\sum_{n=n_0}^\infty e^{-(\frac35|\alpha_n|^2-\delta)x}\Bigr)\le Ce^{-\delta x}.
\end{equation}

For $x\to0^+$  \eqref{expasym} gives us the behavior of $f(x)$. 

For $\Re s>0$ by \eqref{modexp}
\begin{equation}
\Gamma(s/2)Z_\alpha(s)=\int_0^1f(x)x^{\frac{s}{2}}\frac{dx}{x}+
\int_1^\infty f(x)x^{\frac{s}{2}}\frac{dx}{x}:=H_1(s)+H_2(s).
\end{equation}
By \eqref{boundf} $H_2(s)$ is an entire function. 

On the other hand by \eqref{defgN}  
\begin{multline}
H_1(s)=\frac12\int_0^1\Bigl\{\frac{1}{4\sqrt{\pi x}}\log\frac{e^{-C_0}}{16\pi^2 x}+
\frac74
+\sum_{n=1}^{N-1} \frac{a_{n+1}}{\Gamma(1+n/2)}x^{n/2}\Bigr\}x^{\frac{s}{2}}
\frac{dx}{x}+\\
+\frac12\int_0^1 g_N(x)x^{\frac{s}{2}}\frac{dx}{x}
\end{multline}
When $\Re s>1$ we may compute some of the integrals so that
\begin{multline}\label{H1eq}
H_1(s)=\frac{1}{4\sqrt{\pi}}\Bigl\{\frac{2}{(s-1)^2}+\frac{1}{s-1}\log\frac{e^{-C_0}}{16\pi^2}
\Bigr\}+\sum_{n=0}^{N-1} \frac{a_{n+1}}{\Gamma(1+n/2)}\frac{1}{n+s}+\\
+\frac12\int_0^1 g_N(x)x^{\frac{s}{2}}\frac{dx}{x}
\end{multline}
But due to \eqref{ordergN} the last integral defines an analytic function for 
$\Re s>-N$. Hence  \eqref{H1eq} implies that $H_1(s)$ extends to a meromorphic function on $\Re s>N$.

It follows that $\Gamma(s/2)Z_\alpha(s)$ extends to a meromorphic function on  
$\Re s>-N$ for all natural numbers.  Hence $Z_\alpha(s)$ extends to a meromorphic
function on $\C$. 

The function $\Gamma(s/2)Z_\alpha(s)$ has a pole of order $2$ at $s=1$ and poles 
at the points $s=-n$ for $n=0$, $1$, $2$, \dots From \eqref{H1eq} we know the main 
part at  these poles. It follows easily that $Z_\alpha(s)$ has a pole of order $2$ at
$s=1$ and simple poles at the odd negative integers $s=-1$, $s=-3$, \dots 

Finally near $s=1$ we have 
\begin{multline}
Z_\alpha(s)=\frac{1}{\Gamma(s/2)}
\frac{1}{4\sqrt{\pi}}\Bigl\{\frac{2}{(s-1)^2}+\log\frac{e^{-C_0}}{16\pi^2}
\frac{1}{s-1}+\sum_{n=0}^\infty c_n(s-1)^n\Bigr\}=\\
=\Bigl\{\frac{1}{4\pi}+\frac{C_0+\log4}{8\pi}(s-1)+\dots\Bigr\}
\Bigl\{\frac{2}{(s-1)^2}+\log\frac{e^{-C_0}}{16\pi^2}
\frac{1}{s-1}+\dots\Bigr\}=\\
=\frac{1}{2\pi}\Bigl\{\frac{1}{(s-1)^2}-\frac{\log2\pi}{(s-1)}+\cdots\Bigr\}.
\end{multline}
So the function $Z_\alpha(s)$ has a double pole at $s=1$ with residue $-\frac{\log2\pi}{2\pi}$ and main part  given above.

In the same way near $s=-2n$ with $n=0$, $1$, \dots 
\begin{multline}
Z_\alpha(s)=\Bigl\{(-1)^n \,n!	\frac12(s+2n)+\cdots\Bigr\}
\Bigl\{\frac{a_{2n+1}}{\Gamma(1+n)}\frac{1}{2n+s}+\cdots\Bigr\}=\\
(-1)^{n}\frac{a_{2n+1}}{2}+\Orden(s+2n).
\end{multline}
This is equivalent to \eqref{Zvalue}. 

And near the point $s=-(2n-1)$ with $n=1$, $2$, \dots  
\begin{multline}
Z_\alpha(s)=\Bigl\{\frac{1}{\Gamma(-n+\tfrac12)}+\cdots\Bigr\}\Bigl\{
\frac{a_{2n}}{\Gamma(n+\frac12)}\frac{1}{s+2n-1}+\cdots\Bigr\}=\\
=(-1)^{n}\frac{  a_{2n}}{\pi(s+2n-1)}+\cdots\qedhere
\end{multline}
\end{proof}

\end{document}